\documentclass[a4paper,11pt,openright]{article}
\usepackage{}
\usepackage{amsfonts}
\usepackage{mathrsfs}

\usepackage{amssymb}

\newtheorem{lemma}{Lemma}[section]
\newtheorem{theorem}{Theorem}[section]
\newtheorem{corollary}{Corollary}[section]
\newtheorem{definition}{Definition}[section]
\newtheorem{proposition}{Proposition}[section]

\def\blemma{\begin{lemma}\sl{}\def\elemma{\end{lemma}}}
\def\btheorem{\begin{theorem}\sl{}\def\etheorem{\end{theorem}}}
\def\bcorollary{\begin{corollary}\sl{}\def\ecorollary{\end{corollary}}}

\def\bproposition{\begin{proposition}\sl{}\def\eproposition{\end{proposition}}}

\def\beqlb{\begin{eqnarray}}\def\eeqlb{\end{eqnarray}}
\def\beqnn{\begin{eqnarray*}}\def\eeqnn{\end{eqnarray*}}

\def\qed{\hfill$\Box$\medskip}

\def\<{\langle}\def\>{\rangle}

\def\e{{\mbox{\rm e}}}

\begin{document}

\bigskip

\bigskip

{\LARGE\bf Strong law of large number for\\
\indent branching Hunt processes \footnote{The author is supported by NSFC grants (No.
11301020)}}

\bigskip

{\bf Li Wang}\footnote{School of Sciences, Beijing University of
Chemical Technology, Beijing 100029, P.R. China\\
    E-mail\,$:$ wangli@mail.buct.edu.cn }

\bigskip\bigskip

{\narrower{\narrower

\centerline{\bf Abstract}

\bigskip

In this paper we prove that, under certain conditions, a strong law
of large number holds for a class of branching particle systems $X$
corresponding to the parameters $(Y,\beta,\psi)$, where $Y$ is a
Hunt process and $\psi$ is the generating function for the
offspring. The main tool of this paper is the spine decomposition
and we only need a $L\log L$ condition.
\bigskip


\bigskip


\par}\par}

\bigskip\bigskip

\section{Introduction}

\subsection{Motivation}

\setcounter{equation}{0}

In recent years, many people (see
\cite{ET02,EW06,E08,EHK10,CS07,CRW08,LYS13} and the reference
therein) have studied limit theorems for branching Markov processes
or superprocesses using the principal eigenvalue and ground state of
the linear part of the characteristic equations. For superprocesses,
the second moment condition on the branching mechnisms can be weaken,
see \cite{ET02,LYS13}. However, for branching Markov processes, all the
papers in the literature assumed that the branching mechanisms
satisfy a second moment condition or (and), they assume that the
underlying process is symmetric.

In \cite{AH76a}, Asmussen and Hering established a Kesten-Stigum
$LlogL$ type theorem for a class of branching diffusion processes
under a condition which is later called a positive regular property
in \cite{AH76b}. In \cite{LYS09,LYS11}, Liu, Ren and Song
established Kesten-Stigum $LlogL$ type theorem for super-diffusions
and branching Hunt processes respectively. As a natural continuation
of \cite{LYS09}, Liu, Ren and Song give a strong law of large number
for super-diffusions, see \cite{LYS13}. This paper concerns with the
case of branching Markov processes. We establish a strong law of large
numbers for a class of branching Hunt processes. The main tool is
the spine decomposition. We only assume that the branching
mechanisms satisfy a $L\log L$ condition and the underlying
process need not to be symmetric.

We first introduce the setup in this paper.
Let $E$ be a locally compact separable metric space. Denote by
$E_\Delta:=E\cup \{\Delta\}$ the one point compactification of $E$.
Let $\mathcal{B}(E)$ denote both the Borel $\sigma$-fields on $E$
and the space of functions measurable with respect to itself. Write
$\mathcal {B}_b(E)$ (respectively, $\mathcal {B}^+(E))$ for the
space of bounded (respectively, non-negative)
$\mathcal{B}(E)$-measurable functions on $E$. Let $M_p(E)$ be the
space of finite point measures on $E$, that is,
\[
M_p(E)=\bigg\{\sum_{i=1}^n\delta_{x_i}:
n\in\mathbb{N}~\mbox{and}~x_i\in E,i=1,2,\ldots,n\bigg\}.
\]
As usual, $\langle f, \mu\rangle:=\int_E f(x)\mu(dx)$ for any
function $f$ on $E$ and any measure $\mu\in M_p(E)$.

As a continuation of \cite{LYS11}, the model in this paper is the
same as in that paper, we will state it in the next subsection for
reader's convenience.

\subsection{Model}

Let $Y=\{Y_t,\Pi_x,\zeta\}$ be a Hunt process on $E$, where
$\zeta=\inf\{t>0:Y_t=\Delta\}$ is the lifetime of $Y$. Let
$\{P_t,t\geq 0\}$ be the transition semigroup of $Y$:
\[
P_tf(x)=\Pi_x[f(Y_t)]~~~\mbox{for}~f\in\mathcal{B}^+(E).
\]
Let $m$ be a positive Radon measure on $E$ with full support.
$\{P_t, t\geq 0\}$ can be extended to a strongly continuous
semigroup on $L^2(E,m)$. Let $\{\widehat{P}_t, t\geq 0\}$ be the dual
semigroup of $\{P_t, t\geq 0\}$ on $L^2(E,m)$ satisfy
\[
\int_Ef(x)P_tg(x)m(dx)=\int_Eg(x)\widehat{P}_tf(x)m(dx),~~f,g\in
L^2(E,m).
\]
Throught this paper we assume that

\smallskip

\noindent{\bf Assumption 1.1} (i) There exists a family of continuous strictly positive functions
$\{p(t,\cdot,\cdot); t>0\}$ on $E\times E$ such that for any $(t,x)\in (0,\infty)\times E$, we have
\[
P_tf(x)=\int_Ep(t,x,y)f(y)m(dy),~~~\widehat{P}_tf(x)=\int_Ep(t,y,x)f(y)m(dy).
\]
(ii)The semigroups $\{P_t, t\geq 0\}$ and $\{\widehat{P}_t, t\geq 0\}$ are ultracontractive, that is, for
any $t>0$, there exists a constant $c_t>0$ such that
\[
p(t,x,y)\leq c_t~~\mbox{for any}~(x,y)\in E\times E.
\]

Suppose that $\psi\in \mathcal{B}(E\times [-1,1])$ and $\psi$ is the generating function for
each $x\in E$, that is
\[
\psi(x,z)=\sum_{n=0}^\infty p_n(x)z^n,~~|z|\leq 1,
\]
where $p_n(x)\geq 0$ and $\sum_{n=0}^\infty p_n(x)=1$. The branching
system we are going to study determined by the following properties:
\begin{enumerate}
\item The particles in $E$ move independently according to the law of
$Y$, and each particle has a random birth and a random death time.

\item Given the path $Y$ of a particle and given that the particle is alive at time $t$, its probability of dying
in the interval $[t,t+dt)$ is $\beta(Y_t)dt+o(dt)$.
\item When a particle dies at $x\in E$, it splits into $n$ particles at $x$ with probability $p_n(x)$.
The point $\Delta$ is a cemetery. When a particle reaches $\Delta$, it stays at $\Delta$ for ever and there
is no branching at $\Delta$.
\end{enumerate}

We assume that the functions $\beta(x)$ and
$A(x):=\psi'(x,1)=\sum_{n=0}^\infty np_n(x)$ are bounded
$\mathcal{B}(E)$-measurable and that $p_0(x)+p_1(x)=0$ on $E$. The
last condition implies $A(x)\geq 2$ on $E$. The assumption
$p_0(x)=0$ on $E$ is essential for the probabilistic proof of this
paper since we need the spine to be defined for all $t\geq 0$. The
assumption $p_1(x)=0$ on $E$ is just for convenience as the case
$p_1(x)>0$ can be reduced to the case $p_1(x)=0$ by changing the
parameters $\beta$ and $\psi$ of the branching Hunt process.

Let $X_t(B)$ be the number of particles located in $B\in
\mathcal{B}(E)$ at time $t$. Then $X=\{X_t, t\geq 0\}$ is a Markov
process in $M_p(E)$ which is called a $(Y,\beta,\psi)$-branching
process. The process $X$ has probabilities $\{\mathbb{P}_\mu: \mu\in
M_p(E)\}$, and $\mathbb{E}_\mu$ is expectation with respect to
$\mathbb{P}_\mu$. Then we have
 \beqnn
\mathbb{E}_\mu[\<f,X_t\>]=\Pi_\mu[e_{(1-A)\beta}(t)f(Y_t)],~~~f\in
\mathcal{B}_b^+(E), \eeqnn where $e_c(t)=\exp(-\int_0^t c(Y_s)ds)$
for any $c\in\mathcal{B}_b(E)$. We use $\{P_t^{(1-A)\beta},t\geq
0\}$ to denote the following Feynman-Kac semigroup
\[
P_t^{(1-A)\beta}f(x):=\Pi_x[e_{(1-A)\beta}(t)f(Y_t)],~~~f\in
\mathcal{B}(E).
\]
Under Assumption 1.1, we can show that $\{P_t^{(1-A)\beta}\}$ is strongly continuous on $L^2(E,m)$
and for any $t>0$,
$P_t^{(1-A)\beta}$ admits a density $p_t^{(1-A)\beta}(t,x,y)$
which is jointly continuous in $(x,y)$.

Let $\{\widehat{P}_t^{(1-A)\beta},t\geq 0\}$ be the dual semigroup of $\{P_t^{(1-A)\beta},t\geq 0\}$
defined by
\[
\widehat{P}_t^{(1-A)\beta}f(x)=\int_Ep_t^{(1-A)\beta}(t,y,x)f(y)m(dy),~~~f\in \mathcal{B}^+(E).
\]
write $\mathbf{A}$ and $\mathbf{\hat{A}}$ for the generators of
$\{P_t\}$ and $\{\widehat{P}_t\}$. Then the generators of
$\{P_t^{(1-A)\beta}\}$ and $\{\widehat{P}_t^{(1-A)\beta}\}$ can be
formally written as $\mathbf{A}+(A-1)\beta$ and
$\mathbf{\hat{A}}+(A-1)\beta$ respectively.

Let $\sigma(\mathbf{A}+(A-1)\beta)$ and
$\sigma(\mathbf{\hat{A}}+(A-1)\beta)$ be the spectrum of
$\{P_t^{(1-A)\beta}\}$ and $\{\widehat{P}_t^{(1-A)\beta}\}$,
respectively. It follow from Jentzch's Theorem (Theorem V.6.6 on
p.333 of \cite{SH74}) and the strong continuity of
$\{P_t^{(1-A)\beta}\}$ and $\{\widehat{P}_t^{(1-A)\beta}\}$ that the
common value $\lambda_1:= \sup
Re(\sigma(\mathbf{A}+(A-1)\beta))=\sup
Re(\sigma(\mathbf{\hat{A}}+(A-1)\beta))$ is an eigenvalue of
multiplicity $1$ for both $\mathbf{A}+(A-1)\beta$ and
$\mathbf{\hat{A}}+(A-1)\beta$.
Let $\phi$ be an eigenfunction of
$\mathbf{A}+(A-1)\beta$ associated with $\lambda_1$ and $\widetilde{\phi}$
be an eigenfunction of $\mathbf{\hat{A}}+(A-1)\beta$ associated with $\lambda_1$. By
(Proposition 2.3 in \cite{KS08}) we know that $\phi$ and $\widetilde{\phi}$ are strictly
positive and continuous on $E$. We choose $\phi$ and $\widetilde{\phi}$
so that $\int_E\phi\widetilde{\phi}m(dx)=1$. Then
\[
\phi(x)=e^{-\lambda_1t}P_t^{(1-A)\beta}\phi(x),~~\widetilde{\phi}(x)=
e^{-\lambda_1t}\widehat{P}_t^{(1-A)\beta}\widetilde{\phi}(x),~~x\in E.
\]
Throughout this paper we also assume that

\noindent{\bf Assumption 1.2} $\lambda_1>0$ and  $\int_E
\phi^2(y)\widetilde{\phi}(y)m(dy)<\infty$.

The assumption $\lambda_1>0$ is the condition for supercriticality
of the branching Hunt process.

\smallskip

\noindent{\bf Assumption 1.3} The semigroups $\{P_t^{(1-A)\beta}\}$
and $\{\widehat{P}_t^{(1-A)\beta}\}$ are intrinsic ultracontrative,
that is, for any $t>0$ there exists a constant $c_t$ such that
\beqnn p^{(1-A)\beta}(t,x,y)\leq
c_t\phi(x)\widetilde{\phi}(y),~~~x,y\in E. \eeqnn

We refer to \cite{LYS11} for examples satisfy the above assumptions.

\subsection{Spine Decomposition}

For the convenience of state our main result, we shortly recall the
spine decomposition in \cite{LYS11}. First we extend the probability
measure $\mathbb{P}_{\delta_x}$ to a probability measure
$\widetilde{\mathbb{P}}_{\delta_x}$ under which:

\begin{enumerate}
\item a single particle,  $\widetilde{Y}=\{\widetilde{Y}_t\}_{t\geq 0}$, referred to as the spine,
initially starts at $x$ moves according to the measure $\Pi_x$.
\item Given the trajectory $\widetilde{Y}_{\cdot}$, the fission time $\zeta_u$ of node $u$ on the spine
is distributed according to $L^{\beta(\widetilde{Y})}$, where
$L^{\beta(\widetilde{Y})}$ is the law of the Poisson random measure
with intensity $\beta(\widetilde{Y}_t)dt$.
\item At the fission time $\zeta_u$ of node $u$ in the spine, the
single spine particle is replaced by a random number $r_u$ of
offspring with $r_u$ being distributed according to the law
$P(\widetilde{Y}_{\zeta_u})=(p_k(\widetilde{Y}_{\zeta_u}))_{k\geq
1}$.
\item The spine is chosen uniformly from the $r_u$ particles at the fission time of $u$.
\item Each of the remaining $r_u-1$ particles gives rise to independent copys of a $P$-branching Hunt process started at
its space-time point of creation.
\end{enumerate}

Let $\xi=\{\xi_0=\phi, \xi_1, \xi_2,\ldots\}$ be the selected line
of decent in the spine, let $N=(N_t:t\geq 0)$ to denote the counting
process of fission times along the spine. Write $node_t(\xi)$ for
the node in the spine that is alive at time $t$. It is clear that
$node_t(\xi)=\xi_{N_t}$. Define the natural filtration of the motion
and the birth process along the spine by
\[
\mathcal{G}_t:=\sigma((\widetilde{Y}_s, s\leq t), (node_s(\xi):
s\leq t), (\zeta_u, u<\xi_{N_t}), (r_u:u<\xi_{N_t})), \] and define
$\mathcal{G}=\bigcup_{t\geq 0}\mathcal{G}_t$. Let
$\widetilde{\mathcal{F}}_t:=\sigma((X_s,s\leq t), (node_s(\xi):
s\leq t))$ and $\widetilde{\mathcal{F}}=\bigcup_{t\geq
0}\widetilde{\mathcal{F}}_t$. From the spine construction, we know
that
 \beqnn
\mbox{Prob}(u\in \xi)=\prod_{\nu<u}\frac{1}{r_\nu}.
 \eeqnn
It is easy to see that
 \beqlb\label{1.1}
\sum_{u\in L_t}\prod_{\nu<u}\frac{1}{r_\nu}=1.
 \eeqlb
where $L_t$ is the set of particles that are alive at time $t$. For
the definition of $\widetilde{\mathbb{P}}_{\delta_x}$ and the
relations of $\mathbb{P}_{\delta_x}$ with
$\widetilde{\mathbb{P}}_{\delta_x}$, see \cite{LYS11} for details.

Next we define a probability measure
$\widetilde{\mathbb{Q}}_{\delta_x}$ on the branching Hunt process
with a spine. Before that, we need to give some facts concerning
change of measures.

\noindent{\bf Girsanov change of measure} Let $\mathcal {G}_t=\sigma(Y_s; s\leq t)$. Note that
\[
\frac{\phi(Y_t)}{\phi(x)}e^{-\lambda_1t}e_{(1-A)\beta}(t)
\]
is a martingale under $\Pi_x$, and so we can define a martingale change of measure by
\[
\frac{d\Pi_x^\phi}{d\Pi_x}\bigg|_{\mathcal {G}_t}=\frac{\phi(Y_t)}{\phi(x)}e^{-\lambda_1t}e_{(1-A)\beta}(t).
\]
Then $\{Y,\Pi_x^\phi\}$ is a conservative Markov process, and $\phi\widetilde{\phi}$ is a unique invariant
probability measure for the semigroup $\{P_t^{(1-A)\beta}: t\geq 0\}$, that is, for any $f\in \mathcal{B}^+(E)$,
\[
\int_E\phi(x)\widetilde{\phi}(x)P_t^{(1-A)\beta}f(x)m(dx)=\int_Ef(x)\phi(x)\widetilde{\phi}(x)m(dx).
\]
Let $p^\phi(t,x,y)$ be the transition density of $Y$ in $E$ under $\Pi_x^\phi$. Then
\[
p^\phi(t,x,y)=\frac{e^{-\lambda_1t}}{\phi(x)}p^{(1-A)\beta}(t,x,y)\phi(y).
\]
It follows from Theorem 2.8 in \cite{KS08} that, if Assumption 1.3
holds, there exist constants $c>0$ and $\nu>0$ such that
\beqlb\label{1.2} \bigg|\frac{e^{-\lambda_1t}
p^{(1-A)\beta}(t,x,y)}{\phi(x)\widetilde{\phi}(y)}-1\bigg| \leq
ce^{-\nu t}, ~x\in E. \eeqlb which is equivalent to
 \beqlb\label{1.3} \sup_{x\in E}\bigg|\frac{
p^{\phi}(t,x,y)}{\phi(y)\widetilde{\phi}(y)}-1\bigg| \leq ce^{-\nu
t}. \eeqlb

\noindent{\bf Change of measure for Possion process} Suppose that
given a nonnegative measurable function $\beta(Y_t)$, $t\geq 0$, the
Possion process $(n,L^\beta)$ where
$n=\{\{\sigma_i:i=1,2,\ldots,n_t\}:t\geq 0\}$ has instantaneous rate
$\beta(Y_t)$. Further, assume that $n$ is adapted to
$\{\mathcal{L}_t: t\geq 0\}$. Then under the change of measure
\[
\frac{dL^{A\beta}}{dL^\beta}\bigg|_{\mathcal{L}_t}=\prod_{i\leq n_t}A(Y_t)\cdot
\exp\bigg(-\int_0^t((A-1)\beta)(Y_s)ds \bigg)
\]
the process $(n,L^{A\beta})$ is also a Possion process with rate
$A\beta$. See, Chapter 3 in \cite{JS03}.

\medskip

\noindent{\bf The spine construction} Let $\{\mathcal{F}_t: t\geq
0\}$ be the natural filtration generated by $X$. For any $x\in E$,
we define
\[
M_t(\phi)=e^{-\lambda_1t}\frac{\langle X_t, \phi\rangle}{\phi(x)}.
 \]
Then $\{M_t(\phi),t\geq 0\}$ is a nonnegative martingale with respect to $\{\mathcal{F}_t: t\geq 0\}$.
Define the change of measure
\[
\frac{d\widetilde{\mathbb{Q}}_{\delta_x}}{d\widetilde{\mathbb{P}}_{\delta_x}}\bigg|_{\mathcal{F}_t}=M_t(\phi).
\]
Then, under $\widetilde{\mathbb{Q}}_{\delta_x}$, $X$ can be
constructed as follows:
\begin{enumerate}
\item a single particle, $\widetilde{Y}=\{\widetilde{Y}_t\}_{t\geq 0}$, referred to as the spine, initially starts at $x$ moves
according to the measure $\Pi_x^\phi$;
\item The spine undergoes fission into particles at an accelerated intensity $(A\beta)(\widetilde{Y}_t)dt$;
\item At the fission time $\zeta_u$ of node $u$ in the spine, it give birth to $r_u$ particles with size-biased
offspring distribution
$\widehat{P}(\widetilde{Y}_{\zeta_u}):=(\widehat{P}_k(\widetilde{Y}_{\zeta_u}))_{k\geq
1}$, where $\widehat{P}_k(y):=\frac{kp_k(y)}{A(y)}$, $k=1,2,\ldots$,
$y\in E$.
\item The spine is chosen uniformly from the $r_u$ particles at the fission time of $u$.
\item Each of the remaining $r_u-1$ particles gives rise to independent copys of a $P$-branching Hunt process started at
its space-time point of creation.
\end{enumerate}

\btheorem\label{t1.1} (\cite{LYS11}, Theorem 2.9) (Spine
decomposition) We have the following spine decomposition for the
martingale $M_t(\phi)$£º \beqnn
\phi(x)\widetilde{\mathbb{Q}}_{\delta_x}[M_t(\phi)|\mathcal{G}]=e^{-\lambda_1t}\phi(\widetilde{Y}_t)+
\sum_{u<\xi_{N_t}}(r_u-1)\phi(\widetilde{Y}_{\zeta_u})e^{-\lambda_1\zeta_u}
 \eeqnn
\etheorem

Denote by $M_\infty(\phi)$ the almost sure limit of $M_t(\phi)$ as
$t\rightarrow\infty$. In \cite{LYS11}, the author studied the
relationship between the degeneracy property of $M_\infty(\phi)$ and
the function $l$:
 \beqlb\label{1.4} l(x)=\sum_{k=2}^\infty
k\phi(x)\log^+(k\phi(x))p_k(x),~x\in E.
 \eeqlb
 \btheorem\label{t1.2} (\cite{LYS11}, Theorem 1.6) Suppose that $\{X_t: t\geq 0\}$ is
a $(Y,\beta,\psi)$-branching Hunt process and the Assumptions
1.1-1.3 are satisfied. Then $M_\infty(\phi)$ is a non-degenerate
under $\mathbb{P}_\mu$ for any nonzero measure $\mu\in M_p(E)$ if
and only if
 \beqnn
 \int_E\widetilde{\phi}(x)\beta(x)l(x)m(dx)<\infty,
 \eeqnn
where $l$ is defined by (\ref{1.4}).
 \etheorem

\subsection{Main Result}
Define
 \[
W_{t}(\phi):=e^{-\lambda_1t}\< X_t,\phi\>.
 \]
The main goal of this paper is to establish the following almost
sure convergence result.

\btheorem\label{t1.3} Suppose that $\{X_t: t\geq 0\}$ is a
$(Y,\beta,\psi)$-branching Hunt process and the Assumptions 1.1-1.3
are satisfied. If
$\int_E\widetilde{\phi}(x)\beta(x)l(x)m(dx)<\infty$, then there
exists $\Omega_0\subset \Omega$ with full probability (that is,
$\mathbb{P}_{\delta_x}(\Omega_0)=1$ for every $x\in E$) such that,
for every $\omega\in \Omega_0$ and for every bounded Borel
measurable function $f$ on $E$ with compact support whose set of
discontinuous points has zero $m$-measure, we have \beqlb\label{1.5}
\lim_{t\rightarrow\infty}e^{-\lambda_1t}\<
X_t,f\>=W_\infty(\phi)\int_E\widetilde{\phi}(x)f(x)m(dx), \eeqlb
where $W_\infty(\phi)$ is the $\mathbb{P}_{\delta_x}$-almost sure
limit of $e^{-\lambda_1t}\langle X_t,\phi\rangle$. \etheorem

As a consequence of this theorem we immediately get the following

\bcorollary\label{c1.1} (Strong law of large numbers) Suppose that $\{X_t: t\geq 0\}$
is a $(Y,\beta,\psi)$-branching Hunt process and the Assumptions
1.1-1.3 are satisfied. If $\int_E\widetilde{\phi}(x)\beta(x)l(x)m(dx)<\infty$, then there exists
$\Omega_0\subset \Omega$ with full probability such that, for every $\omega\in \Omega_0$
and for every relatively
compact Borel subset $B$ in $E$ having $m(B)>0$ and $m(\partial B)=0$, we have
$\mathbb{P}_{\delta_x}$-almost surely,
\[\lim_{t\rightarrow\infty}\frac{X_t(B)(\omega)}{\mathbb{P}_{\delta_x}[X_t(B)]}
=\frac{W_\infty(\phi)(\omega)}{\phi(x)}.\] \ecorollary

 \section{Proof of Theorem 1.3}

 \setcounter{equation}{0}

We will prove the theorem by the following steps.

\bproposition\label{p2.1} If
$\int_E\widetilde{\phi}(x)\beta(x)l(x)m(dx)<\infty,$ then for any
$m\in\mathbb{N},\sigma>0$,
\[
\lim_{n\rightarrow\infty}|U_{(m+n)\sigma}(\phi
f)-\mathbb{E}_{\delta_x}(U_{(m+n)\sigma}(\phi f)|\mathcal
{F}_{n\sigma})|=0,~\mathbb{P}_{\delta_x}\mbox{-a.s.}
\]
where $U_{t}(f\phi):=e^{-\lambda_1t}\< X_t,f\phi\>$ for $f\in
\mathcal{B}_b^+(E)$. \eproposition

We will prove this result later. According to the spine
construction, if a particle $u\in \xi$,  then at the fission time
$\zeta_u$, it give birth to $r_u$ offspring, one of which continues
the spine while the other $r_u-1$ individuals go off to create
subtrees which are copies of the original branching Hunt process, we
write them by $(\tau,M)_j^u$, $j=1,\ldots,r_u-1$. Put
\[
X_{t-\zeta_u}^j=\sum_{\nu\in L_t,
\nu\in(\tau,M)_j^u}\delta_{Y_\nu(t)}(\cdot),~t\geq\zeta_u,
\]
where $\{Y_\nu: \nu\in(\tau,M)_j^u\}$ are independent copies of $Y$.
$(X_{t-\zeta_u}^j,t\geq\zeta_u)$ is a $(Y,\beta,\psi)$-branching Hunt process
with birth time $\zeta_u$ and starting point
$\widetilde{Y}_{\zeta_u}$. Then we can write
 \beqnn
U_{t}(f\phi)=e^{-\lambda_1t}(f\phi)(\widetilde{Y}_t)+
e^{-\lambda_1t}\sum_{u<\xi_{N_t}}\sum_{j=1}^{r_u-1}\<f\phi,X_{t-\zeta_u}^j\>.
 \eeqnn Define
 \beqnn \widetilde{U}_{t}(f\phi)=e^{-\lambda_1t}(f\phi)(\widetilde{Y}_t)+
e^{-\lambda_1t}\sum_{u<\xi_{N_t}}\sum_{j=1}^{r_u-1}\<f\phi,X_{t-\zeta_u}^j\>
I_{\{r_u\phi(\widetilde{Y}_{\zeta_u})\leq e^{\lambda_1\zeta_u}\}}.
 \eeqnn
and \[M_t^{u,j}(\phi):=e^{-\lambda_1(t-\zeta_u)}
\frac{\<\phi,X_{t-\zeta_u}^j\>}{\phi(\widetilde{Y}_{\zeta_u})},~~t\geq
\zeta_u. \] Then $\{M_t^{u,j}(\phi),t\geq \zeta_u\}$ is, conditional
on $\mathcal {G}$, a nonnegative
$\widetilde{\mathbb{P}}_{\delta_x}$-martingale on the subtree
$(\tau,M)_j^u$, and therefore
 \beqlb\label{2.1}
\widetilde{\mathbb{Q}}_{\delta_x}[M_t^{u,j}(\phi)|\mathcal{G}]=
\widetilde{\mathbb{P}}_{\delta_x}[M_t^{u,j}(\phi)|\mathcal{G}]=1.
\eeqlb

Suppose that $\{Y_i:i=1,\ldots,L_{n\sigma}\}$ describes the path of
particles alive at time $n\sigma$. Note that we may always write
\beqnn U_{(m+n)\sigma}(f\phi)=\sum_{i=1}^{L_{n\sigma}}
e^{-\lambda_1n\sigma}U^{(i)}_{m\sigma}(f\phi)
 \eeqnn
where given $\mathcal{F}_{n\sigma}$, the collection
$\{U^{(i)}_{m\sigma}(f\phi):i=1,\ldots,L_{n\sigma}\}$ are mutually
independent and equal in distribution to $U_{m\sigma}(f\phi)$ under
$\mathbb{P}_{\delta_{Y_i}}$. Then we can write
 \beqlb\label{2.2} \nonumber
U_{(m+n)\sigma}(f\phi)&=&\sum_{i=1}^{L_{n\sigma}}
e^{-\lambda_1n\sigma}\widetilde{U}^{(i)}_{m\sigma}(f\phi)+\sum_{i=1}^{L_{n\sigma}}
e^{-\lambda_1n\sigma}\left(U^{(i)}_{m\sigma}(f\phi)-\widetilde{U}^{(i)}_{m\sigma}(f\phi)\right)
\\&:=&U^{[1]}_{(m+n)\sigma}(f\phi)+U^{[2]}_{(m+n)\sigma}(f\phi),
\eeqlb
where $U^{[1]}_{(m+n)\sigma}(f\phi)$ and $U^{[2]}_{(m+n)\sigma}(f\phi)$ stand for the first term and the
second term on the right hand respectively.

\blemma\label{l2.1} If
$\int_E\widetilde{\phi}(x)\beta(x)l(x)m(dx)<\infty$ and $\int_E
\phi^2(y)\widetilde{\phi}(y)m(dy)<\infty$. Then,
 for $f\in
\mathcal{B}_b^+(E)$ and $x\in E$, \beqnn
\widetilde{\mathbb{E}}_{\delta_x}[\widetilde{U}_t(\phi f)]^2<\infty
\eeqnn \elemma

\noindent{\it Proof.} First, we rewrite $\widetilde{U}_{t}(f\phi)$
into a new form and take the conditional expectation,
 \beqnn & &\widetilde{\mathbb{E}}_{\delta_x}[\widetilde{U}_{t}(\phi
f)|\mathcal{F}_t]\\
&=&\sum_{u\in L_t}\left(e^{-\lambda_1t}(f\phi)(\widetilde{Y}_u(t)) +
e^{-\lambda_1t}\sum_{\nu<u}\sum_{j=1}^{r_\nu-1}\<f\phi,X_{t-\zeta_\nu}^j\>
I_{\{r_\nu\phi(\widetilde{Y}_{\zeta_\nu})\leq
e^{\lambda_1\zeta_\nu}\}}\right)\widetilde{\mathbb{E}}_{\delta_x}(I_{\{u\in
\xi\}}|\mathcal{F}_t)\\
&=&\sum_{u\in L_t}\left(e^{-\lambda_1t}(f\phi)(\widetilde{Y}_u(t)) +
e^{-\lambda_1t}\sum_{\nu<u}\sum_{j=1}^{r_\nu-1}\<f\phi,X_{t-\zeta_\nu}^j\>
I_{\{r_\nu\phi(\widetilde{Y}_{\zeta_\nu})\leq
e^{\lambda_1\zeta_\nu}\}}\right)\prod_{\nu<u}\frac{1}{r_\nu}\\
&\stackrel{d}{=}&e^{-\lambda_1t}(f\phi)(\widetilde{Y}_t)+
e^{-\lambda_1t}\sum_{\nu<
\xi_{N_t}}\sum_{j=1}^{r_\nu-1}\<f\phi,X_{t-\zeta_\nu}^j\>
I_{\{r_\nu\phi(\widetilde{Y}_{\zeta_\nu})\leq
\e^{\lambda_1\zeta_\nu}\}},
 \eeqnn
 where in the last equation, in order not to introduce another
symbol, we still use $\xi_{N_t}$ to denote one of the particles
alive at time $t$, ``$\stackrel{d}{=}$" means equal in distribution
under $\mathbb{P}_{\delta_x}$ and the equality (\ref{1.1}) was used.
Using (\ref{2.1}) and $\widetilde{U}_{t}(\phi f)\leq
\|f\|_\infty\cdot W_t(\phi)$, we have
 \beqnn
&&\phi(x)^{-1}\widetilde{\mathbb{E}}_{\delta_x}[\widetilde{U}_{t}(\phi
f)]^2\\&\leq &\|f\|_\infty\cdot
\phi(x)^{-1}\widetilde{\mathbb{E}}_{\delta_x}[W_{t}(\phi)\widetilde{U}_{t}(\phi
f)]\\
&=&\|f\|_\infty\cdot
\phi(x)^{-1}\widetilde{\mathbb{E}}_{\delta_x}\left[W_{t}(\phi)\widetilde{\mathbb{E}}_{\delta_x}[\widetilde{U}_{t}(\phi
f)|\mathcal{F}_t]\right]\\
&=&\|f\|_\infty\cdot\widetilde{\mathbb{Q}}_{\delta_x}\left[e^{-\lambda_1t}(f\phi)(\widetilde{Y}_t)+
e^{-\lambda_1t}\sum_{\nu<
\xi_{N_t}}\sum_{j=1}^{r_\nu-1}\<f\phi,X_{t-\zeta_\nu}^j\>
I_{\{r_\nu\phi(\widetilde{Y}_{\zeta_\nu})\leq
e^{\lambda_1\zeta_\nu}\}}\right]\\
&\leq&\|f\|_\infty^2\cdot\widetilde{\mathbb{Q}}_{\delta_x}\bigg[e^{-\lambda_1t}\phi(\widetilde{Y}_t)+\sum_{
\nu<\xi_{N_t}}(r_\nu-1)\phi(\widetilde{Y}_{\zeta_\nu})e^{-\lambda_1\zeta_\nu}
I_{\{r_\nu\phi(\widetilde{Y}_{\zeta_\nu})\leq e^{\lambda_1\zeta_\nu}\}}\bigg]\\
&\leq&
\|f\|_\infty^2\left(\Pi_x^\phi[e^{-\lambda_1t}\phi(\widetilde{Y}_t)]+\Pi_x^\phi\int_0^t
e^{-\lambda_1s}\phi(\widetilde{Y}_s)\beta(\widetilde{Y}_s)A(\widetilde{Y}_s)\sum_{k=2}^\infty
k \widehat{p}_k(\widetilde{Y}_s)I_{\{k\phi(\widetilde{Y}_s)\leq
e^{\lambda_1 s}\}} ds\right).
 \eeqnn
where $\|\cdot\|_\infty$ means the supremum norm here and in the
paper. Call the two expressions in bracket
 on the right hand side the \textit{spine
term} $A(x,t)$ and the \textit{sum term} $B(x,t)$ respectively. Note
that (\ref{1.3}) implies that \beqnn \bigg|\int_E
p^{\phi}(t,x,y)\phi(y)m(dy)-\int_E
\phi^2(y)\widetilde{\phi}(y)m(dy)\bigg|\leq ce^{-\nu t}\int_E
\phi^2(y)\widetilde{\phi}(y)m(dy). \eeqnn Therefore,
 \beqlb\label{2.3}
e^{\lambda_1t}A(x,t)=\Pi_x^\phi (\phi(\widetilde{Y}_t))\leq (c+1)\int_E
\phi^2(y)\widetilde{\phi}(y)m(dy)<\infty.
 \eeqlb
For the sum term, using the assumption that $A$ and $\beta$ are
bounded, we get
 \beqlb\label{2.4}
\nonumber B(x,t)&=&\Pi_x^\phi\int_0^t
e^{-\lambda_1s}\phi(\widetilde{Y}_s)\beta(\widetilde{Y}_s)A(\widetilde{Y}_s)\sum_{k=2}^\infty
k \widehat{p}_k(\widetilde{Y}_s)I_{\{k\phi(\widetilde{Y}_s)\leq e^{\lambda_1
s}\}}ds\\ \nonumber &=&\Pi_x^\phi\int_0^t
e^{-\lambda_1s}\phi(\widetilde{Y}_s)\beta(\widetilde{Y}_s)A(\widetilde{Y}_s)\sum_{k=2}^\infty
k\frac{k}{A(\widetilde{Y}_s)}p_k(\widetilde{Y}_s)I_{\{k\phi(\widetilde{Y}_s)\leq
e^{\lambda_1 s}\}}ds\\ \nonumber &\leq&\Pi_x^\phi\int_0^t
e^{-(\lambda_1-\lambda_1)s}\beta(\widetilde{Y}_s)\sum_{k=2}^\infty k
p_k(\widetilde{Y}_s)I_{\{k\phi(\widetilde{Y}_s)\leq e^{\lambda_1 s}\}}
ds\\  &\leq& \|\beta A\|_\infty\cdot t<\infty,
 \eeqlb
for all $x\in E$, then the conclusion follows. \qed

\blemma\label{l2.2} If
$\int_E\widetilde{\phi}(x)\beta(x)l(x)m(dx)<\infty,$ then for any
$m\in\mathbb{N},\sigma>0$,
 \beqnn
 & &\sum_{n=1}^\infty\mathbb{\widetilde{P}}_{\delta_x}
 \{U_{(n+m)\sigma}(f\phi)\neq U^{[1]}_{(n+m)\sigma}(f\phi)\}<\infty,\\
&
&\sum_{n=1}^\infty\widetilde{\mathbb{E}}_{\delta_x}\left[U^{[1]}_{(m+n)\sigma}(\phi
f)-\widetilde{\mathbb{E}}_{\delta_x}(U^{[1]}_{(m+n)\sigma}(\phi
f)|\mathcal {\widetilde{F}}_{n\sigma})\right]^2<\infty.
 \eeqnn
 where $U^{[1]}_{(m+n)\sigma}(f\phi)$ was defined in (\ref{2.2}). In particular
\[
\lim_{n\rightarrow\infty}\left|U^{[1]}_{(m+n)\sigma}(\phi
f)-\mathbb{\widetilde{E}}_{\delta_x}(U^{[1]}_{(m+n)\sigma}(\phi
f)|\mathcal{\widetilde{F}}_{n\sigma})\right|=0,~\mathbb{\widetilde{P}}_{\delta_x}\mbox{-a.s.}
\]
\elemma {\it Proof.} Note that (\ref{1.2}) implies that for any
$s\in [0, m\sigma]$, there is a constant $C_{m\sigma}$ such that
\beqnn p^{(1-A)\beta}(s,x,y)\leq
C_{m\sigma}\phi(x)\widetilde{\phi}(y),~~~x,y\in E. \eeqnn Then by
the spine construction and Fubini theorem, we get \beqnn
&&\sum_{n=1}^\infty\mathbb{\widetilde{P}}_{\delta_x}
 \{U_{(n+m)\sigma}(f\phi)\neq U^{[1]}_{(n+m)\sigma}(f\phi)\}\\
&\leq&\sum_{n=1}^\infty \mathbb{\widetilde{E}}_{\delta_x}\left[
 \sum_{i=1}^{L_{n\sigma}}
\mathbb{\widetilde{P}}_{\delta_x}\left(U^{(i)}_{m\sigma}(f\phi)
\neq\widetilde{U}^{(i)}_{m\sigma}(f\phi)|\mathcal
{\widetilde{F}}_{n\sigma}\right)\right]\\
&\leq& \sum_{n=1}^\infty \mathbb{\widetilde{E}}_{\delta_x}\left(
 \sum_{i=1}^{L_{n\sigma}}
\int_0^{m\sigma}ds\int_Ep^{(1-A)\beta}(s,Y_i,y)\beta(y)\sum_{k=2}^\infty
p_k(y)I_{\{k\phi(y)>e^{\lambda_1 (s+n\sigma)}\}}m(dy) \right)\\
&\leq&C_{m\sigma}\sum_{n=1}^\infty
\mathbb{\widetilde{E}}_{\delta_x}\left(
 \sum_{i=1}^{L_{n\sigma}}
\int_0^{m\sigma}ds\int_E\phi(Y_i)\widetilde{\phi}(y)\beta(y)\sum_{k=2}^\infty
p_k(y)I_{\{k\phi(y)>e^{\lambda_1 (s+n\sigma)}\}}m(dy) \right)\\
&=&C_{m\sigma}\sum_{n=1}^\infty
e^{\lambda_1n\sigma}\phi(x)
\left(\int_0^{m\sigma}ds\int_E\widetilde{\phi}(y)\beta(y)\sum_{k=2}^\infty
p_k(y)I_{\{k\phi(y)>e^{\lambda_1 (s+n\sigma)}\}}m(dy)\right)\\
&=&C_{m\sigma}\phi(x)\left(
\int_0^{m\sigma}ds\int_E\widetilde{\phi}(y)\beta(y)\sum_{n=1}^\infty
e^{\lambda_1n\sigma}\sum_{k=2}^\infty
p_k(y)I_{\{k\phi(y)>e^{\lambda_1 (s+n\sigma)}\}}m(dy)\right)\\
&\leq&C_{m\sigma}\phi(x)\left(
\int_0^{m\sigma}ds\int_E\widetilde{\phi}(y)\beta(y)\sum_{k=2}^\infty
\sum_{n=1}^{\frac{1}{\lambda_1\sigma}log^+[k\phi(y)]}
e^{\lambda_1n\sigma} p_k(y)I_{\{k\phi(y)>e^{\lambda_1
(s+n\sigma)}\}}m(dy)\right)\\
&\leq&\frac{C_{m\sigma}}{\lambda_1\sigma}\phi(x)\left(
\int_0^{m\sigma}ds\int_E\widetilde{\phi}(y)\beta(y)\sum_{k=2}^\infty
log^+[k\phi(y)]k\phi(y)p_k(y)m(dy)\right)\\
&=&\frac{C_{m\sigma}m}{\lambda_1}\phi(x)
\int_E\widetilde{\phi}(y)\beta(y)l(y)m(dy)<\infty. \eeqnn

\noindent For the second inequality, recall that, if $X_i$ are
independent random variables with $E(X_i)=0$ or they are martingale
difference, then \beqnn E\bigg|\sum_{i=1}^n X_i\bigg|^p\leq
2^p\sum_{i=1}^n E|X_i|^p. \eeqnn Jensen's inequality also implies
that $|u+v|^p\leq 2^{p-1}(|u|^p+|v|^p)$ for $p\in (1,2]$. Then we
have \beqnn & &\mathbb{E}(|U_{s+t}-\mathbb{E}(U_{s+t}|\mathcal
{F}_{t})|^p|\mathcal
{F}_{t}))\\
&\leq&
2^pe^{-\lambda_1pt}\sum_{i=1}^{L_t}\mathbb{E}\bigg(|U^{(i)}_{s}-\mathbb{E}(U^{(i)}_{s}|\mathcal
{F}_{t})|^p\bigg|\mathcal
{F}_{t}\bigg)\\
&\leq&
2^pe^{-\lambda_1pt}\sum_{i=1}^{L_t}\mathbb{E}\bigg(2^{p-1}(|U^{(i)}_{s}|^p+|\mathbb{E}(U^{(i)}_{s}|\mathcal
{F}_{t})|^p)\bigg|\mathcal
{F}_{t}\bigg)\\
&\leq&
2^pe^{-\lambda_1pt}\sum_{i=1}^{L_t}2^{p-1}\mathbb{E}\bigg(|U^{(i)}_{s}|^p+|\mathbb{E}(U^{(i)}_{s}|\mathcal
{F}_{t})|^p)\bigg|\mathcal
{F}_{t}\bigg)\\
&\leq&
2^{2p}e^{-\lambda_1pt}\sum_{i=1}^{L_t}\mathbb{E}(|U^{(i)}_{s}|^p|\mathcal
{F}_{t}). \eeqnn Note that for any $f\in \mathcal{B}_b^+(E)$,
$U_t(f\phi)\leq \|f\|_\infty\cdot W_t(\phi)$, we have that \beqnn &
&\sum_{n\geq
1}\widetilde{\mathbb{E}}_{\delta_x}\left(\left|U^{[1]}_{(m+n)\sigma}-\mathbb{E}
\left(\widetilde{U}^{[1]}_{(m+n)\sigma}\bigg|\widetilde{\mathcal
{F}}_{n\sigma}\right)\right|^2\right)\\
&\leq& 2^4\sum_{n\geq
1}e^{-2\lambda_1n\sigma}\widetilde{\mathbb{E}}_{\delta_x}\bigg(\sum_{i=1}^{L_{n\sigma}}
\widetilde{\mathbb{E}}_{\delta_{Y_i}}[\widetilde{U}^{(i)}_{m\sigma}(\phi
f)]^2\bigg)\\
&\leq& 2^4\sum_{n\geq
1}\mathbb{E}_{\delta_x}\bigg(\sum_{i=1}^{L_{n\sigma}}e^{-2\lambda_1n\sigma}
\phi(Y_i)(A(Y_i,m\sigma)+B(Y_i,m\sigma))\bigg)\\
&=&2^4\sum_{n\geq 1}e^{-\lambda_1n\sigma}\phi(x)\Pi_x^\phi
[A(Y_{n\sigma},m\sigma)+B(Y_{n\sigma},m\sigma)]
 \eeqnn
where $A(x,t)$ and $B(x,t)$ were defined in Lemma \ref{l2.1}. Then
as a consequence of the previous estimates (\ref{2.3}) and
(\ref{2.4}), we conclude that the last sum remains finite.\qed

\blemma\label{l2.3} If
$\int_E\widetilde{\phi}(x)\beta(x)l(x)m(dx)<\infty,$ then for any
$m\in\mathbb{N},\sigma>0$,
 \beqnn
 \sum_{n=0}^\infty\widetilde{\mathbb{E}}_{\delta_x}\left[\left(U_{(m+n)\sigma}(\phi
f)-U^{[1]}_{(m+n)\sigma}(\phi f)\right)\bigg|\mathcal
{\widetilde{F}}_{n\sigma}\right]
~\mbox{converges}~\mathbb{\widetilde{P}}_{\delta_x}\mbox{-a.s.}
 \eeqnn
\elemma {\it Proof.} Take $f=1$ in Lemma \ref{l2.2}, then $\{U_t(\phi): t\geq 0\}$ is a nonnegative
martingale. By Lemma \ref{l2.2} we have
 \beqlb
 & &\label{2.5}\sum_{n=0}^\infty\mathbb{\widetilde{P}}_{\delta_x}
 \left\{U_{(n+m)\sigma}(\phi)\neq U^{[1]}_{(n+m)\sigma}(\phi)\right\}<\infty,\\
&
&\label{2.6}\sum_{n=0}^\infty\widetilde{\mathbb{E}}_{\delta_x}\left[U^{[1]}_{(m+n)\sigma}(\phi
)-\widetilde{\mathbb{E}}_{\delta_x}\left(U^{[1]}_{(m+n)\sigma}(\phi
)\bigg|\mathcal {\widetilde{F}}_{n\sigma}\right)\right]^2<\infty.
 \eeqlb
Note that
 \beqnn
\widetilde{\mathbb{E}}_{\delta_x}\left[U^{[1]}_{(m+n)\sigma}(\phi
)\bigg|\mathcal{\widetilde{F}}_{n\sigma}\right]&=&
\widetilde{\mathbb{E}}_{\delta_x}\left[\left(U_{(m+n)\sigma}(\phi
)-U^{[2]}_{(m+n)\sigma}(\phi
)\right)\bigg|\mathcal{\widetilde{F}}_{n\sigma}\right]\\
&=&U_{n\sigma}(\phi)-\widetilde{\mathbb{E}}_{\delta_x}\left[U^{[2]}_{(m+n)\sigma}(\phi
)\bigg|\mathcal {\widetilde{F}}_{n\sigma}\right]
 \eeqnn
By (\ref{2.5}) and (\ref{2.6}), we have
 \beqnn
\sum_{n=0}^\infty \left(U_{(m+n)\sigma}(\phi
)-U_{n\sigma}(\phi)+\widetilde{\mathbb{E}}_{\delta_x}\left[U^{[2]}_{(m+n)\sigma}(\phi
)\bigg|\mathcal
{\widetilde{F}}_{n\sigma}\right]\right)\mbox{converges}~\mathbb{\widetilde{P}}_{\delta_x}\mbox{-a.s.}
 \eeqnn
since $U_t(\phi)$ converges almost surely as $t\rightarrow\infty$,
we have
 \beqnn
\sum_{n=0}^\infty\widetilde{\mathbb{E}}_{\delta_x}\left[U^{[2]}_{(m+n)\sigma}(\phi
)\bigg|\mathcal {\widetilde{F}}_{n\sigma}\right]
~\mbox{converges}~\mathbb{\widetilde{P}}_{\delta_x}\mbox{-a.s.}
 \eeqnn
So we have
 \beqnn
\sum_{n=0}^\infty\widetilde{\mathbb{E}}_{\delta_x}\left[\left(U_{(m+n)\sigma}(\phi
f)-U^{[1]}_{(m+n)\sigma}(\phi f)\right)\bigg|\mathcal
{\widetilde{F}}_{n\sigma}\right] \leq\|f\|_\infty
\sum_{n=0}^\infty\widetilde{\mathbb{E}}_{\delta_x}\left[U^{[2]}_{(m+n)\sigma}(\phi
)\bigg|\mathcal {\widetilde{F}}_{n\sigma})\right] \eeqnn converges
$\mathbb{\widetilde{P}}_{\delta_x}$-a.s. \qed

\medskip

\noindent {\it Proof of Porposition \ref{p2.1}.} From the
decomposition (\ref{2.2}), we have
 \beqnn
&&U_{(m+n)\sigma}(\phi f)-\mathbb{E}_{\delta_x}(U_{(m+n)\sigma}(\phi
f)|\mathcal
{F}_{n\sigma})\\
&=& U_{(m+n)\sigma}(\phi
f)-\widetilde{\mathbb{E}}_{\delta_x}(U_{(m+n)\sigma}(\phi
f)|\widetilde{\mathcal {F}}_{n\sigma})\\
&=&U_{(m+n)\sigma}(\phi f)-U^{[1]}_{(m+n)\sigma}(\phi f)
+U^{[1]}_{(m+n)\sigma}(\phi
f)-\widetilde{\mathbb{E}}_{\delta_x}\left(U^{[1]}_{(m+n)\sigma}(\phi
f)\bigg|\mathcal {\widetilde{F}}_{n\sigma}\right)\\
&
&-~\widetilde{\mathbb{E}}_{\delta_x}\left[\left(U_{(m+n)\sigma}(\phi
f)-U^{[1]}_{(m+n)\sigma}(\phi f)\right)\bigg|\mathcal
{\widetilde{F}}_{n\sigma}\right] \eeqnn Now the conclusion of this
proposition follows immediately form Lemma \ref{l2.1}, Lemma\ref{l2.2} and
Lemma \ref{l2.3}.\qed

\btheorem\label{t2.1} If
$\int_E\widetilde{\phi}(x)\beta(x)l(x)m(dx)<\infty$, then
for any $\sigma>0$ and $f\in \mathcal{B}_b^+(E)$,
\[
\lim_{n\rightarrow\infty}e^{-\lambda_1n\sigma}\<\phi
f,X_{n\sigma}\>=W_\infty(\phi)\int_E\widetilde{\phi}(x)\phi(x)f(x)m(dx),~\mathbb{P}_{\delta_x}\mbox{-a.s.}
\]
\etheorem {\it Proof.} By Markov property of branching processes we
have
 \beqnn
\mathbb{E}_\mu[e^{-\lambda_1(m+n)\sigma}\<\phi f,
X_{(m+n)\sigma}\>|\mathcal
{F}_{n\sigma}]=e^{-\lambda_1n\sigma}\<e^{-\lambda_1m\sigma}
P^{(1-A)\beta}_{m\sigma}(\phi f),X_{n\sigma}\>. \eeqnn Note that
(\ref{1.2}) implies that, for any $f\in \mathcal{B}_b^+(E)$, \beqnn
\bigg|\frac{e^{-\lambda_1m\sigma} P^{(1-A)\beta}_{m\sigma}(\phi
f)(x)}{\phi(x)}-\int_E \phi(y)\widetilde{\phi}(y)
f(y)m(dy)\bigg|\leq ce^{-\nu m\sigma}\int_E
\phi(y)\widetilde{\phi}(y) f(y)m(dy), \eeqnn which is equivalent to
\beqnn \bigg|\frac{e^{-\lambda_1m\sigma}
P^{(1-A)\beta}_{m\sigma}(\phi f)(x)}{\phi(x)\int_E
\phi(y)\widetilde{\phi}(y) f(y)m(dy)}-1\bigg|\leq ce^{-\nu m\sigma}.
\eeqnn Thus there exist positive constants $c_m\leq 1$ and $C_m\geq
1$ such that \beqnn
 c_m\phi(x)\int_E \phi(y)\widetilde{\phi}(y)
f(y)m(dy)\leq e^{-\lambda_1m\sigma} P^{(1-A)\beta}_{m\sigma}(\phi f)(x)\leq C_m\phi(x)\int_E \phi(y)\widetilde{\phi}(y)
f(y)m(dy),
\eeqnn
and that $\lim_{m\rightarrow\infty}c_m=\lim_{m\rightarrow\infty}C_m=1$. Hence,
\beqnn
e^{-\lambda_1n\sigma}\<e^{-\lambda_1m\sigma}
P^{(1-A)\beta}_{m\sigma}(\phi f),X_{n\sigma}\>
&\geq&c_me^{-\lambda_1n\sigma}\<\phi,X_{n\sigma}\>\int_E \phi(y)\widetilde{\phi}(y)f(y)m(dy)\\
&=& c_m W_{n\sigma}(\phi)\int_E \phi(y)\widetilde{\phi}(y)
f(y)m(dy),
\eeqnn
and
\beqnn
e^{-\lambda_1n\sigma}\<e^{-\lambda_1m\sigma}
P^{(1-A)\beta}_{m\sigma}(\phi f),X_{n\sigma}\>
&\leq&C_me^{-\lambda_1n\sigma}\<\phi,X_{n\sigma}\>\int_E \phi(y)\widetilde{\phi}(y)f(y)m(dy)\\
&=& C_m W_{n\sigma}(\phi)\int_E \phi(y)\widetilde{\phi}(y)
f(y)m(dy). \eeqnn Those two inequalities and Proposition \ref{p2.1}
imply that \beqnn
\limsup_{n\rightarrow\infty}e^{-\lambda_1n\sigma}\<\phi
f,X_{n\sigma}\>
&=&\limsup_{n\rightarrow\infty}e^{-\lambda_1(m+n)\sigma}\<\phi f, X_{(m+n)\sigma}\>\\
&=&\limsup_{n\rightarrow\infty}e^{-\lambda_1n\sigma}\<e^{-\lambda_1m\sigma}
P^{(1-A)\beta}_{m\sigma}(\phi f),X_{n\sigma}\>\\
&\leq&\limsup_{n\rightarrow\infty}C_m W_{n\sigma}(\phi)\int_E \phi(y)\widetilde{\phi}(y)
f(y)m(dy)\\
&=&C_m W_\infty(\phi)\int_E \phi(y)\widetilde{\phi}(y)
f(y)m(dy),~~\mathbb{P}_{\delta_x}\mbox{-a.s.}
\eeqnn
and that
\beqnn
\liminf_{n\rightarrow\infty}e^{-\lambda_1n\sigma}\<\phi f,X_{n\sigma}\>\geq c_m W_\infty(\phi)\int_E \phi(y)\widetilde{\phi}(y)
f(y)m(dy),~~\mathbb{P}_{\delta_x}\mbox{-a.s.}
\eeqnn
Letting $m\rightarrow\infty$, we get
\beqnn
\lim_{n\rightarrow\infty}e^{-\lambda_1n\sigma}\<\phi f,X_{n\sigma}\>=W_\infty(\phi)\int_E \phi(y)\widetilde{\phi}(y)
f(y)m(dy),~~\mathbb{P}_{\delta_x}\mbox{-a.s.}
\eeqnn
The proof is now complete. \qed

 \blemma\label{l2.4} If
$\int_E\widetilde{\phi}(x)\beta(x)l(x)m(dx)<\infty$, then
 for any open subset $U$ in $E$ and $x\in E$, we have
\beqnn \liminf_{t\rightarrow\infty}e^{-\lambda_1 t}\<\phi
I_U,X_t\>\geq W_\infty(\phi)\int_E \phi(y)\widetilde{\phi}(y)I_U(y)
m(dy),~~\mathbb{P}_{\delta_x}\mbox{-a.s.} \eeqnn \elemma {\it
Proof.} For $x\in E$ and $\varepsilon>0$, let
\[
U^\varepsilon(x):=\bigg\{y\in U: \phi(y)\geq
\frac{1}{1+\varepsilon}\phi(x)\bigg\}.
\]
Define
\[
Z_{n,\nu}^{\sigma,\varepsilon}=\frac{1}{1+\varepsilon}~\phi(Y_\nu(n\sigma))1_{\{Y_\nu(t)\in
U^\varepsilon(Y_\nu(n\sigma)) ~\mbox{for every}~t\in
[n\sigma,(n+1)\sigma)\}}
\]
where each $Y_\nu$ describes the motion of particle $\nu$ in the
branching particle system. Let
 \[
S_n^{\sigma,\varepsilon}=e^{-\lambda_1
n\sigma}\sum_{u<\xi_{N_{n\sigma}}}\sum_{j=1}^{r_u-1}\sum_{\nu\in(\tau,M)_j^u}Z_{n,\nu}^{\sigma,\varepsilon}
 \]
where the subtrees $\{(\tau,M)_j^u\}$ were defined below Proposition
\ref{p2.1}.
For $t\in [n\sigma,(n+1)\sigma)$, we have
 \beqnn
e^{-\lambda_1 t}\<\phi I_U,X_t\>&=&e^{-\lambda_1 t}(\phi
I_U)(\widetilde{Y}_t)+e^{-\lambda_1
t}\sum_{u<\xi_{N_t}}\sum_{j=1}^{r_u-1}\sum_{\nu\in(\tau,M)_j^u}\phi(Y_\nu(t))1_{\{Y_\nu(t)\in
U\}}\\
&\geq&e^{-\lambda_1 t}(\phi
I_U)(\widetilde{Y}_t)+\frac{e^{-\lambda_1
n\sigma}}{1+\varepsilon}\sum_{u<\xi_{N_{n\sigma}}}\sum_{j=1}^{r_u-1}
\sum_{\nu\in(\tau,M)_j^u}e^{-\lambda_1\sigma}\phi(Y_\nu(n\sigma))1_{\{Y_\nu(t)\in
U^\varepsilon(Y_\nu(n\sigma))\}}.
 \eeqnn
Applying Proposition \ref{p2.1} with $
\<f\phi,X_{(m+n)\sigma-\zeta_u}^j\>\left(=\sum_{\nu\in(\tau,M)_j^u}(f\phi)(Y_\nu((n+m)\sigma))\right)
$
 replaced by
\[\sum_{\nu\in(\tau,M)_j^u}e^{-\lambda_1\sigma}\phi(Y_\nu(n\sigma))1_{\{Y_\nu(t)\in
U^\varepsilon(Y_\nu(n\sigma)),~\mbox{for
every}~t\in[n\sigma,(n+1)\sigma)\}},
\]
similar estimates to those found in Proposition \ref{p2.1} show us
that
\[
\lim_{n\rightarrow\infty}|S_n^{\sigma,\varepsilon}-E(S_n^{\sigma,\varepsilon}|\mathcal
{F}_{n\sigma})|=0,~\mathbb{P}_{\delta_x}\mbox{-a.s.}
\]
Then we have
\beqnn
\liminf_{t\rightarrow\infty}e^{-\lambda_1 t}\<\phi I_U,X_t\>&\geq&
e^{-\lambda_1 \sigma}\liminf_{n\rightarrow\infty}S_n^{\sigma,\varepsilon}\\
&=&e^{-\lambda_1\sigma}\liminf_{n\rightarrow\infty}e^{-\lambda_1
n\sigma}\sum_{i=1}^{L_{n\sigma}}
E_{Y^i_{n\sigma}}[S_0^{\sigma,\varepsilon}]\\
&=&e^{-\lambda_1\sigma}W_\infty(\phi)\int_E
\widetilde{\phi}(y)\xi_U^{\sigma,\varepsilon}(y) m(dy),
\eeqnn where we have used Theorem \ref{t2.1} in the last equality
and \beqnn
\xi_U^{\sigma,\varepsilon}(Y^i_{n\sigma})=E_{Y^i_{n\sigma}}
\left[S_0^{\sigma,\varepsilon}\right]=\frac{\phi(Y^i_{n\sigma})}
{1+\varepsilon}P_{Y^i_{n\sigma}}\bigg(Y(t)\in
U^\varepsilon(Y^i_{n\sigma})~\mbox{for all}~t\in[0,\sigma)\bigg).
\eeqnn Taking $\sigma\downarrow 0$, we get that $\int_E
\widetilde{\phi}(y)\xi_U^{\sigma,\varepsilon}(y)
m(dy)\rightarrow \frac{1}{1+\varepsilon}\int_E
\phi(y)\widetilde{\phi}(y)1_U(y) m(dy)$; hence subsequently taking
$\varepsilon\downarrow 0$ gives us \beqlb\label{2.7}
\liminf_{t\rightarrow\infty}e^{-\lambda_1 t}\<\phi I_U,X_t\>\geq
W_\infty(\phi)\int_E \phi(y)\widetilde{\phi}(y)1_U(y) m(dy).
\eeqlb\qed

\noindent{\bf Proof of Theorem \ref{t1.3}} Since $E$ is a locally
compact separable metric space, there exists a countable base
$\mathcal{U}$ of open set $\{U_k,k\geq 1\}$ that is closed under
finite union. By Lemma \ref{l2.4}, there exists $\Omega_0\subset
\Omega$ of full probability so that for every $\omega\in\Omega_0$,
\beqnn \liminf_{t\rightarrow\infty}e^{-\lambda_1 t}\<\phi
I_{U_k},X_t\>\geq W_\infty(\phi)\int_{U_k}
\phi(y)\widetilde{\phi}(y) m(dy). \eeqnn For any open set $U$, there
exists a sequence of increasing open sets $\{U_{n_k},k\geq 1\}$ in
$\mathcal{U}$ so that $\bigcup_{k=1}^\infty U_{n_k}=U$. We have for
every $\omega\in\Omega_0$, \beqnn
\liminf_{t\rightarrow\infty}e^{-\lambda_1 t}\<\phi I_U,X_t\>&\geq&
\liminf_{t\rightarrow\infty}
e^{-\lambda_1 t}\<\phi I_{U_{n_k}},X_t\>\\
&\geq& W_\infty(\phi)\int_{U_{n_k}}
\phi(y)\widetilde{\phi}(y)m(dy)~\mbox{for every}~k\geq 1. \eeqnn
Passing $k\rightarrow\infty$ yields that \beqnn
\liminf_{t\rightarrow\infty}e^{-\lambda_1 t}\<\phi I_{U},X_t\>\geq
W_\infty(\phi)\int_{U} \phi(y)\widetilde{\phi}(y) m(dy). \eeqnn We
consider (\ref{1.5}) on $\{W_\infty(\phi)>0\}$. For each fixed $\omega\in
\Omega_0\cap \{W_\infty(\phi)>0\}$, define the probability measure $\mu_t$
and $\mu$ on $E$ respectively, by
\[\mu_t(A)(\omega)=\frac{e^{-\lambda_1 t}\<\phi I_{A},X_t\>}{W_t(\phi)(\omega)}~~\mbox{and}~~\mu(A)=\int_A\phi(x)\widetilde{\phi}(x)
m(dx),~~A\in\mathcal{B}(E)
\]
for every $t\geq 0$. Note that the measure $\mu_t$ is well defined
for every $t\geq 0$. The inequality (\ref{2.7}) tell us that $\mu_t$
converges weakly to $\mu$. Since $\phi$ is strictly positive and
continuous on $E$, for every function $f$ on $E$ with compact
support whose discontinuity set has zero $m$-measure (equivalently
zero $\mu$-measure), $h:=\frac{f}{\phi}$ is a bounded function
having compact support with the same set of discontinuity with $f$.
We thus have
\[
\lim_{t\rightarrow\infty}\int_E hd\mu_t=\int_E hd\mu
\]
which is equivalent to say that
\beqnn
\lim_{t\rightarrow\infty}e^{-\lambda_1t}\< X_t,f\>(\omega)=W_\infty(\phi)(\omega)
\int_E\widetilde{\phi}(x)f(x)m(dx)~~~~
\mbox{for every}~\omega\in \Omega_0\cap \{W_\infty(\phi)>0\}.
\eeqnn
Since, for every function $f$ on $E$ such that $|f|$ is bounded by $c\phi$ for some $c>0$,
\[
e^{-\lambda_1t}|\< X_t,f\>|\leq e^{-\lambda_1t}\< X_t,|f|\>\leq cW_t,
\]
(\ref{1.5}) holds automatically on $\{W_\infty(\phi)=0\}$. This complete
the proof of the theorem.\qed

\noindent{\bf Proof of corollary \ref{c1.1}} It is enough if we can prove that
 \beqnn
 \lim_{t\rightarrow\infty}e^{-\lambda_1t}\mathbb{E}_{\delta_x}\<
X_t,f\>=\phi(x)\int_E\widetilde{\phi}(x)f(x)m(dx),~~~f\in \mathcal{B}_b^+(E).
 \eeqnn
Note that, (\ref{1.3}) implies that $\Pi_x^\phi (f(Y_t))\rightarrow
\<f\phi,\widetilde{\phi}\>$ for every measurable function $f$
satisfying $\<f\phi,\widetilde{\phi}\><\infty$. Then we have \beqnn
e^{-\lambda_1t}\mathbb{E}_{\delta_x}\<
X_t,f\>&=&e^{-\lambda_1t}P_t^{(1-A)\beta}f(x)\\
&=&e^{-\lambda_1t}\int_Ep_t^{(1-A)\beta}(t,x,y)f(y)m(dy)\\
&=&\phi(x)\int_Ep_t^{\phi}(t,x,y)\frac{f(y)}{\phi(y)}m(dy)\\
&\rightarrow&\phi(x)\int_E\widetilde{\phi}(x)f(x)m(dx).
\eeqnn
Combining with Theorem \ref{t1.3}, we get the desired result.\qed

\bigskip

\textbf{Acknowledgement}. I would like to give my sincere thanks to
 my supervisor Professor Zenghu Li for his encouragement and helpful
 discussions. Thanks are also given to the referee for her or his
 valuable comments and suggestions.

\end{document}